\documentclass[11pt,reqno]{amsart}
   \usepackage[centertags]{amsmath}
   \usepackage{amsfonts}
   \usepackage{amsmath}
   \usepackage{amssymb}
   \usepackage{amsthm}
   \usepackage{newlfont}
   \usepackage[latin1]{inputenc}
\usepackage{multirow, bigdelim}

\usepackage{epsfig}
\usepackage{pst-all}
\usepackage{pstricks}
\usepackage{pgf}
\usepackage{mathrsfs}
\usepackage{hyperref}
\usepackage{bm}

\usepackage{graphicx}

\usepackage{bm}

\usepackage{mathtools}
\usepackage{ifmtarg}

\usepackage{lscape}

\makeatletter

\allowdisplaybreaks[3]

\theoremstyle{plain}

\theoremstyle{definition}

\theoremstyle{remark}

\newtheorem*{remark*}{Remark}

\numberwithin{equation}{section}

\title[An urn model for the Jacobi-Pi\~neiro polynomials]{An urn model for the Jacobi-Pi\~neiro polynomials}

\author{F. Alberto Gr\"unbaum}
\address{F. Alberto Gr\"unbaum\\
Department of Mathematics. University of California, Berkeley. Berkeley, CA 94720  U.S.A.}
\email{grunbaum@math.berkeley.edu}
\date{\today}
\thanks{
}

\author{Manuel D. de la Iglesia}
\address{Manuel D. de la Iglesia\\
Instituto de Matem\'aticas, Universidad Nacional Aut\'onoma de M\'exico, Circuito Exterior, C.U., 04510, Ciudad de M\'exico, M\'exico.}
\email{mdi29@im.unam.mx}
\thanks{This work of the second author was partially supported by PAPIIT-DGAPA-UNAM grant IN104219 (M\'exico) and CONACYT grant A1-S-16202 (M\'exico).}

\date{\today}



\begin{document}

\maketitle

\begin{abstract}
The list of physically motivated urn models that can be solved in terms of classical orthogonal polynomials is very small. It includes a model proposed by D. Bernoulli and further analyzed by S. Laplace and a model proposed by P. and T. Ehrenfest and eventually connected with the Krawtchouk and Hahn polynomials. This connection was reversed recently in the case of the Jacobi polynomials where a rather contrived, and later a simpler urn model was proposed. Here we consider an urn model going with the Jacobi-Pi\~neiro multiple orthogonal polynomials. These polynomials have recently been put forth in connection with a stochastic matrix.
\end{abstract}

\section{Introduction}

Urn models for the flow of two incompressible liquids between two containers or for the exchange of heat between two isolated bodies, which are nowadays given as examples of finite state Markov chains, predate any formal introduction of these probabilistic notions and even a rigorous set-up for probability theory as given by A. Kolmogorov around 1930. Two celebrated examples are the Bernoulli-Laplace (1769-1812) and the Ehrenfest (1907, see \cite{E}) models discussed in W. Feller's book (see \cite{ F}, pages 121 and 378).

At a much later time it was noticed that these physically motivated Markov chains can be analyzed completely in terms of families of orthogonal polynomials that occupy important places in the Askey-Wilson tableaux, namely Hahn and Krawtchouk polynomials. For an exposition of this see \cite{G2} and its references, or more recently the monograph \cite{MDIB}.

The situation described above was part of the motivation for \cite{G4} where an urn model (albeit an elaborate one) was proposed for the Markov chain with an infinite number of states going with the Jacobi polynomials. In this case one is dealing with a discrete-time birth-death chain stemming from a semi-infinite tridiagonal matrix. This approach was adapted to the case of matrix-valued Jacobi polynomials (see \cite{GdI4, GPT3}) where one is dealing with a quasi-birth-and-death process (see \cite{LaR, DRSZ}). This results in a stochastic matrix that is a semi-infinite block tridiagonal matrix. The importance of studying different walks by means of ``spectral methods'' has been exploited beyond classical walks. In the case of quantum walks, results such as recurrence are related to spectral properties of a measure on the unit circle, see \cite{GVWW,BGVW}.

In \cite{GdI3} a new idea was introduced and exploited: by using a stochastic LU (or UL) factorization of the stochastic matrix going with the Jacobi polynomials a much simpler model than the one in \cite{G4} was put forward. This idea has been adapted to other situations such as multivariate orthogonal polynomials (see \cite{FdI}).

The purpose of this paper is to use the approach in \cite{GdI3} to give an urn model associated to the Jacobi-Pi\~neiro polynomials. They originate in \cite{Pin} and give a basic example of so-called multiple orthogonal polynomials. For references on this kind of polynomials the reader can consult \cite{Pin, Nik, vAC, ABvA}. These are polynomials in one variable indexed by two non-negative indices $(n_1,n_2)$. They give rise to a family of polynomials indexed by one non-negative index $n$ (see Proposition 9 of \cite{ManP}) and these polynomials satisfy a higher-order recursion relation. The fact that this gives rise to a stochastic matrix has been recently pointed out in \cite{ManP}.

\section{Stochastic LU factorization of the stochastic matrix}

Let $\{X_t: t=0,1,\ldots\}$ be the Markov chain on $\mathbb{Z}_{\geq0}$ with transition probability matrix given by the stochastic matrix $P_{II}$ which appears in Corollary 7 of \cite{ManP}. Here we will denote by $P$ this matrix, which is associated with the type II Jacobi-Pi\~neiro multiple orthogonal polynomials. We also denote by $a_n,b_n,c_{n+1},d_{n+2},n\geq0,$ the corresponding transition probabilities given by
\begin{equation}\label{tpm}
\begin{split}
a_n&=\mathbb{P}\left(X_{t+1}=n+1\, |\, X_t=n\right),\quad n\geq0,\\
b_n&=\mathbb{P}\left(X_{t+1}=n \, |\, X_t=n\right),\quad n\geq0,\\
c_n&=\mathbb{P}\left(X_{t+1}=n-1 \, |\, X_t=n\right),\quad n\geq1,\\
d_n&=\mathbb{P}\left(X_{t+1}=n-2 \, |\, X_t=n\right),\quad n\geq2.
\end{split}
\end{equation}
Therefore $P$ is the semi-infinite banded matrix given by 
\begin{equation}
P=\begin{pmatrix}\label{PP}
b_0&a_0&0&&\\
c_1&b_1&a_1&0&&\\
d_2&c_2&b_2&a_2&0&\\
0&d_3&c_3&b_3&a_3&\\
&&\ddots&\ddots&\ddots&\ddots
\end{pmatrix}.
\end{equation}
A diagram of the transitions between the states is given by

\vspace{1.4cm}

\begin{center}
$$\begin{psmatrix}[colsep=1.9cm]
  \cnode{.4}{0}& \cnode{.4}{1} & \cnode{.4}{2}& \cnode{.4}{3}& \cnode{.4}{4}& \cnode{.4}{5} &  \rnode{6}{\Huge{\cdots}} \\
\psset{nodesep=3pt,arcangle=15,labelsep=2ex,linewidth=0.3mm,arrows=->,arrowsize=1mm
3} \nccurve[angleA=160,angleB=200,ncurv=4]{0}{0}
\uput[u](0.95,1.95){a_0}\uput[u](2.85,1.95){a_1}\uput[u](4.75,1.95){a_2}\uput[u](6.65,1.95){a_3}\uput[u](8.55,1.95){a_4}\uput[u](10.65,1.95){a_5}
\uput[u](0.95,1.1){c_1}\uput[u](2.85,1.1){c_2}\uput[u](4.75,1.1){c_3}\uput[u](6.65,1.1){c_4}\uput[u](8.55,1.1){c_5}\uput[u](10.65,1.1){c_6}
\uput[u](-0.75,1.95){b_0}\uput[u](1.9,2.7){b_1}
\uput[u](3.8,0.25){b_2}
\uput[u](5.7,2.7){b_3}\uput[u](7.6,0.25){b_4}\uput[u](9.5,2.7){b_5}
\uput[u](1.9,3.35){d_2}\uput[u](5.7,3.35){d_4}\uput[u](9.5,3.35){d_6}\uput[u](3.8,-0.4){d_3}\uput[u](7.6,-0.4){d_5}
\nccurve[angleA=70,angleB=110,ncurv=4]{1}{1}
\nccurve[angleA=-70,angleB=-110,ncurv=4]{2}{2}
\nccurve[angleA=70,angleB=110,ncurv=4]{3}{3}
\nccurve[angleA=-70,angleB=-110,ncurv=4]{4}{4}
\nccurve[angleA=70,angleB=110,ncurv=4]{5}{5}
\nccurve[angleA=115,angleB=65,ncurv=1]{2}{0}
\nccurve[angleA=115,angleB=65,ncurv=1]{4}{2}
\nccurve[angleA=115,angleB=65,ncurv=1]{6}{4}
\nccurve[angleA=-115,angleB=-65,ncurv=1]{3}{1}
\nccurve[angleA=-115,angleB=-65,ncurv=1]{5}{3}
 \ncarc{0}{1}\ncarc{1}{0} \ncarc{1}{2} \ncarc{2}{1}\ncarc{2}{3} \ncarc{3}{2}
\ncarc{3}{4} \ncarc{4}{3} \ncarc{4}{5}\ncarc{5}{4} \ncarc{5}{6} \ncarc{6}{5}
\psset{labelsep=-3.40ex}\nput{90}{0}{0}
\psset{labelsep=-3.40ex}\nput{90}{1}{1}
\psset{labelsep=-3.40ex}\nput{90}{2}{2}
\psset{labelsep=-3.40ex}\nput{90}{3}{3}
\psset{labelsep=-3.40ex}\nput{90}{4}{4}
\psset{labelsep=-3.40ex}\nput{90}{5}{5}
\end{psmatrix}
$$
\end{center}
\vspace{.1cm}

Following \cite{ManP}, $P$ depends on three parameters $\alpha,\beta$ and $\gamma$, and it is stochastic if and only if $\alpha,\beta,\gamma>-1, \alpha\neq\beta$ and $|\alpha-\beta|<1$ (although the case $\alpha=\beta$ is allowed too, as we will see below). The matrix $P$ \eqref{PP} admits a \emph{stochastic} LU factorization, as in \cite{GdI3}, of the form
\begin{equation}\label{Pxy}
P=P_LP_U=\begin{pmatrix}
s_0&0&\\
r_1&s_1&0&\\
t_2&r_2&s_2&0&\\
0&t_3&r_3&s_3&0&\\
&&\ddots&\ddots&\ddots&\ddots
\end{pmatrix}\begin{pmatrix}
y_0&x_0&0&\\
0&y_1&x_1&0&\\
&0&y_2&x_3&0&\\
&&&\ddots&\ddots&\ddots
\end{pmatrix},
\end{equation}
where, for $n\geq0$, we have
\begin{equation}\label{xy}
\begin{split}
x_{2n}&=\frac{2n+\gamma+1}{3n+\alpha+\gamma+2},\quad y_{2n}=\frac{n+\alpha+1}{3n+\alpha+\gamma+2},\\
x_{2n+1}&=\frac{2n+\gamma+2}{3n+\beta+\gamma+3},\quad y_{2n+1}=\frac{n+\beta+1}{3n+\beta+\gamma+3},
\end{split}
\end{equation}
\begin{equation}\label{ts}
\begin{split}
t_{2n}&=\frac{n(n+\alpha-\beta)}{(3n+\alpha+\gamma)(3n+\alpha+\gamma+1)},\; s_{2n}=\frac{(2n+\alpha+\gamma+1)(2n+\beta+\gamma+1)}{(3n+\alpha+\gamma+1)(3n+\beta+\gamma+1)},\\
t_{2n+1}&=\frac{n(n+\beta-\alpha)}{(3n+\beta+\gamma+1)(3n+\beta+\gamma+2)},\; s_{2n+1}=\frac{(2n+\alpha+\gamma+2)(2n+\beta+\gamma+2)}{(3n+\alpha+\gamma+3)(3n+\beta+\gamma+2)},
\end{split}
\end{equation}
and
\begin{equation}\label{rr}
\begin{split}
r_{2n}&=\frac{n(2n+\beta+\gamma)}{(3n+\alpha+\gamma)(3n+\alpha+\gamma+1)}+\frac{n(2n+\alpha+\gamma+1)}{(3n+\alpha+\gamma+1)(3n+\beta+\gamma+1)},\\
r_{2n+1}&=\frac{n(2n+\alpha+\gamma+1)}{(3n+\beta+\gamma+1)(3n+\beta+\gamma+2)}+\frac{(n+1)(2n+\beta+\gamma+2)}{(3n+\alpha+\gamma+3)(3n+\beta+\gamma+2)}.
\end{split}
\end{equation}
Notice that all the coefficients are nonnegative, $x_n+y_n=1$ and $t_n+s_n+r_n=1$ for $n\geq0$, i.e. the matrices $P_L$ and $P_U$ are stochastic. Observe that the factor $P_U$ is a pure-birth  chain on $\mathbb{Z}_{\geq0}$ with diagram
\begin{center}
\vspace{0.8cm}
$$
\begin{psmatrix}[colsep=1.9cm]
  \cnode{.4}{0}& \cnode{.4}{1} & \cnode{.4}{2}& \cnode{.4}{3}& \cnode{.4}{4}& \cnode{.4}{5} &  \rnode{6}{\Huge{\cdots}} \\
\psset{nodesep=3pt,arcangle=15,labelsep=2ex,linewidth=0.3mm,arrows=->,arrowsize=1mm
3} \nccurve[angleA=160,angleB=200,ncurv=4]{0}{0}
\uput[u](0.95,1.95){x_0}\uput[u](2.85,1.95){x_1}\uput[u](4.75,1.95){x_2}\uput[u](6.65,1.95){x_3}\uput[u](8.55,1.95){x_4}\uput[u](10.65,1.95){x_5}
\uput[u](-0.75,1.95){y_0}\uput[u](1.9,2.7){y_1}\uput[u](3.8,2.7){y_2}\uput[u](5.7,2.7){y_3}\uput[u](7.6,2.7){y_4}\uput[u](9.5,2.7){y_5}
\nccurve[angleA=70,angleB=110,ncurv=4]{1}{1}
\nccurve[angleA=70,angleB=110,ncurv=4]{2}{2}
\nccurve[angleA=70,angleB=110,ncurv=4]{3}{3}
\nccurve[angleA=70,angleB=110,ncurv=4]{4}{4}
\nccurve[angleA=70,angleB=110,ncurv=4]{5}{5}
 \ncarc{0}{1}\ncarc{1}{2}\ncarc{2}{3}
\ncarc{3}{4}\ncarc{4}{5} \ncarc{5}{6}
\psset{labelsep=-3.40ex}\nput{90}{0}{0}
\psset{labelsep=-3.40ex}\nput{90}{1}{1}
\psset{labelsep=-3.40ex}\nput{90}{2}{2}
\psset{labelsep=-3.40ex}\nput{90}{3}{3}
\psset{labelsep=-3.40ex}\nput{90}{4}{4}
\psset{labelsep=-3.40ex}\nput{90}{5}{5}
\end{psmatrix}
$$
\end{center}
\vspace{-1.0cm}
while $P_L$ is a pure-death chain on $\mathbb{Z}_{\geq0}$ with absorbing state at 0 with diagram
\vspace{1.8cm}
$$\begin{psmatrix}[colsep=1.9cm]
  \cnode{.4}{0}& \cnode{.4}{1} & \cnode{.4}{2}& \cnode{.4}{3}& \cnode{.4}{4}& \cnode{.4}{5} &  \rnode{6}{\Huge{\cdots}} \\
\psset{nodesep=3pt,arcangle=15,labelsep=2ex,linewidth=0.3mm,arrows=->,arrowsize=1mm
3} \nccurve[angleA=160,angleB=200,ncurv=4]{0}{0}
\uput[u](0.95,1.1){r_1}\uput[u](2.85,1.1){r_2}\uput[u](4.75,1.1){r_3}\uput[u](6.65,1.1){r_4}\uput[u](8.55,1.1){r_5}\uput[u](10.65,1.1){r_6}
\uput[u](-0.75,1.95){1}\uput[u](1.9,2.7){s_1}\uput[u](5.7,2.7){s_3}\uput[u](9.5,2.7){s_5}
\uput[u](3.8,0.25){s_2}
\uput[u](7.6,0.25){s_4}
\uput[u](1.9,3.35){t_2}\uput[u](5.7,3.35){t_4}\uput[u](9.5,3.35){t_6}\uput[u](3.8,-0.4){t_3}\uput[u](7.6,-0.4){t_5}
\nccurve[angleA=70,angleB=110,ncurv=4]{1}{1}
\nccurve[angleA=-70,angleB=-110,ncurv=4]{2}{2}
\nccurve[angleA=70,angleB=110,ncurv=4]{3}{3}
\nccurve[angleA=-70,angleB=-110,ncurv=4]{4}{4}
\nccurve[angleA=70,angleB=110,ncurv=4]{5}{5}
\nccurve[angleA=115,angleB=65,ncurv=1]{2}{0}
\nccurve[angleA=115,angleB=65,ncurv=1]{4}{2}
\nccurve[angleA=115,angleB=65,ncurv=1]{6}{4}
\nccurve[angleA=-115,angleB=-65,ncurv=1]{3}{1}
\nccurve[angleA=-115,angleB=-65,ncurv=1]{5}{3}
\ncarc{1}{0} \ncarc{2}{1}\ncarc{3}{2}
\ncarc{4}{3} \ncarc{5}{4} \ncarc{6}{5}
\psset{labelsep=-3.40ex}\nput{90}{0}{0}
\psset{labelsep=-3.40ex}\nput{90}{1}{1}
\psset{labelsep=-3.40ex}\nput{90}{2}{2}
\psset{labelsep=-3.40ex}\nput{90}{3}{3}
\psset{labelsep=-3.40ex}\nput{90}{4}{4}
\psset{labelsep=-3.40ex}\nput{90}{5}{5}
\end{psmatrix}
$$
\vspace{0.2cm}
The factorization \eqref{Pxy} gives the following relations among all coefficients
\begin{equation}\label{abcd}
\begin{split}
a_n&=s_{n}x_n,\quad n\geq0\\
b_n&= r_{n}x_{n-1}+ s_ny_n,\quad n\geq1,\quad b_0=s_0y_0,\\
c_n&=r_n y_{n-1}+t_nx_{n-2},\quad n\geq2,\quad c_1=r_1y_0,\\
c_n&=t_ny_{n-2},\quad n\geq2,
\end{split}
\end{equation}
which is an alternative way of writing the transition probabilities appearing in Corollary 7 of \cite{ManP}. Now we can see from the coefficients in \eqref{xy}--\eqref{rr} that the case $\alpha=\beta$  can be allowed (although it may be exceptional in terms of the Jacobi-Pi\~neiro polynomials).

In \cite{ManP} one finds a different matrix $P_I$ going with the type I Jacobi-Pi\~neiro polynomials. Its shape is given as that of the transpose of the one in \eqref{PP}. The methods and models used here for $P_{II}$ can be easily extended to the matrix $P_I$.

\section{The urn model}

In this section we assume $\alpha=1/M,\beta=1/N,\gamma\geq0, M,N,\gamma\in\mathbb{Z}_{>0}$. The condition $|\alpha-\beta|<1$ is now equivalent to $|M-N|<MN$. After this substitution we have that the coefficients in \eqref{xy}, \eqref{ts} and \eqref{rr} are given, for $n\geq0$, by
\begin{equation}\label{xy2}
\begin{split}
x_{2n}&=\frac{M(2n+\gamma+1)}{3Mn+M\gamma+2M+1},\quad y_{2n}=\frac{M(n+1)+1}{3Mn+M\gamma+2M+1},\\
x_{2n+1}&=\frac{N(2n+\gamma+2)}{3Nn+N\gamma+3N+1},\quad y_{2n+1}=\frac{N(n+1)+1}{3Nn+N\gamma+3N+1},
\end{split}
\end{equation}

\begin{equation}\label{ts2}
\begin{split}
t_{2n}&=\frac{Mn(MNn+N-M)}{N(3Mn+M\gamma+1)(3Mn+M\gamma+M+1)},\\
s_{2n}&=\frac{(2Mn+M\gamma+M+1)(2Nn+N\gamma+N+1)}{(3Mn+M\gamma+M+1)(3Nn+N\gamma+N+1)},\\
t_{2n+1}&=\frac{Nn(MNn+M-N)}{M(3Nn+N\gamma+N+1)(3Nn+N\gamma+2N+1)},\\
s_{2n+1}&=\frac{(2Mn+M\gamma+2M+1)(2Nn+N\gamma+2N+1)}{(3Mn+M\gamma+3M+1)(3Nn+N\gamma+2N+1)},
\end{split}
\end{equation}
and
\begin{equation}\label{rr2}
\begin{split}
r_{2n}&=\frac{M^2n(2Nn+N\gamma+1)}{N(3Mn+M\gamma+1)(3Mn+M\gamma+M+1)}\\
&\qquad\qquad+\frac{Nn(2Mn+M\gamma+M+1)}{(3Mn+M\gamma+M+1)(3Nn+N\gamma+N+1)},\\
r_{2n+1}&=\frac{N^2n(2Mn+M\gamma+M+1)}{M(3Nn+N\gamma+N+1)(3Nn+N\gamma+2N+1)}\\
&\qquad\qquad+\frac{M(n+1)(2Nn+N\gamma+2N+1)}{(3Mn+M\gamma+3M+1)(3Nn+N\gamma+2N+1)}.
\end{split}
\end{equation}

Consider the LU factorization $P=P_LP_U$ \eqref{Pxy} where the entries of $P_L$ and $P_U$ are given in \eqref{xy2}, \eqref{ts2} and \eqref{rr2}. Each one of these matrices $P_L$ and $P_U$ will give rise to an experiment in terms of an urn model, which we call Experiment 1 and Experiment 2, respectively. At times $t = 0, 1, 2,\ldots$ an urn A contains $m$ blue balls and this determines the state of our Markov chain at that time. Additionally, for Experiment 1, we will have two other urns, one painted in blue, which we call B, and the other one painted in red, which we call R. These auxiliary urns are initially empty and will be emptied after their use in going from one time step to the next. All the urns sit in a bath consisting of an infinite number of blue and red balls. These urns are not used in Experiment 2.

Let us consider first the (slighty simpler) Experiment 2 (for $P_U$), and let us call this chain $\{X_t^{(2)}, t=0,1,\ldots\}$. Assume that urn A contains $m$ blue balls ($m\geq0$) at time 0 (i.e. $X_0^{(2)}=m$). Consider first the case where $m$ is even and write $m=2n, n\geq0$. Remove all the balls from urn A and place in A $M(2n+\gamma+1)$ blue balls and $M(n+1)+1$ red balls from the bath. Draw one ball from urn A at random with the uniform distribution. We have two possibilities:
\begin{itemize}
\item If we get a blue ball then we remove into the bath or add from the bath balls until we have $m+1=2n+1$ blue balls in urn A and start over. Therefore
$$
\mathbb{P}\left(X_1^{(2)}=2n+1 \, | \,  X_0^{(2)}=2n\right)=\frac{M(2n+\gamma+1)}{3Mn+M\gamma+2M+1}.
$$
Observe that this probability is given by $x_{2n}$ in \eqref{xy2}.

\item If we get a red ball then we remove into the bath or add from the bath balls until we have $m=2n$ blue balls in urn A and start over. Therefore
$$
\mathbb{P}\left(X_1^{(2)}=2n \, | \,  X_0^{(2)}=2n\right)=\frac{M(n+1)+1}{3Mn+M\gamma+2M+1}.
$$
Observe that this probability is given by $y_{2n}$ in \eqref{xy2}.
\end{itemize}

Consider now the case where $m$ is odd and write $m=2n+1, n\geq0$. Remove all the balls from urn A and place in A $N(2n+\gamma+2)$ blue balls and $N(n+1)+1$ red balls from the bath. Draw one ball from urn A at random with the uniform distribution. We have two possibilities:
\begin{itemize}
\item If we get a blue ball then we remove into the bath or add from the bath balls until we have $m+1=2n+2$ blue balls in urn A and start over. Therefore
$$
\mathbb{P}\left(X_1^{(2)}=2n+2 \, | \,  X_0^{(2)}=2n+1\right)=\frac{N(2n+\gamma+2)}{3Nn+N\gamma+3N+1}.
$$
Observe that this probability is given by $x_{2n+1}$ in \eqref{xy2}.

\item If we get a red ball then we remove into the bath or add from the bath balls until we have $m=2n+1$ blue balls in urn A and start over. Therefore
$$
\mathbb{P}\left(X_1^{(2)}=2n+1 \, | \,  X_0^{(2)}=2n+1\right)=\frac{N(n+1)+1}{3Nn+N\gamma+3N+1}.
$$
Observe that this probability is given by $y_{2n+1}$ in \eqref{xy2}.
\end{itemize}

Now we describe Experiment 1 (for $P_L$), and let us call this chain $\{X_t^{(1)}, t=0,1,\ldots\}$. Recall that now we will use two auxiliary urns B and R which are initially empty and will be emptied after their use in each time step. Assume that urn A contains $m$ blue balls ($m\geq0$) at time $0$. Consider first the case where $m$ is even and write $m=2n,n\geq1$ (the state 0 is absorbing). Remove all balls from urn A and place in A $Mn$ blue balls and $2Mn+M\gamma+M+1$ red balls from the bath. In urn B we place $MNn+N-M$ blue balls and $M(2Nn+N\gamma+1)$ red balls. In urn R we place $Nn$ blue balls and $2Nn+N\gamma+N+1$ red balls. Draw one ball from urn A at random with the uniform distribution. If we get a blue ball then we go to urn B and draw a ball, while if we get a red ball then we go urn R and draw a ball. We have three possibilities:
\begin{itemize}
\item If we get two blue balls in a row, i.e. one from urn A and then one from urn B, then we remove into the bath or add from the bath balls until we have $m-2=2n-2$ blue balls in urn A and start over. Therefore
$$
\mathbb{P}\left(X_1^{(1)}=2n-2 \, | \,  X_0^{(1)}=2n\right)=\frac{Mn}{3Mn+M\gamma+M+1}\cdot\frac{MNn+N-M}{N(3Mn+M\gamma+1)}.
$$
Observe that this probability is given by $t_{2n}$ in \eqref{ts2}.
\item  If we get two red balls in a row, i.e. one from urn A and then one from urn R, then we remove into the bath or add from the bath balls until we have $m=2n$ blue balls in urn A and start over. Therefore
$$
\mathbb{P}\left(X_1^{(1)}=2n \, | \,  X_0^{(1)}=2n\right)=\frac{2Mn+M\gamma+M+1}{3Mn+M\gamma+M+1}\cdot\frac{2Nn+N\gamma+N+1}{3Nn+N\gamma+N+1}.
$$
Observe that this probability is given by $s_{2n}$ in \eqref{ts2}.
\item If we get either a blue and a red ball or a red and a blue ball then we remove into the bath or add from the bath balls until we have $m-1=2n-1$ blue balls in urn A and start over. Therefore
\begin{align*}
\mathbb{P}\left(X_1^{(1)}=2n-1 \, | \,  X_0^{(1)}=2n\right)&=\frac{Mn}{3Mn+M\gamma+M+1}\cdot\frac{M(2Nn+N\gamma+1)}{N(3Mn+M\gamma+1)}\\
&\qquad\qquad+\frac{2Mn+M\gamma+M+1}{3Mn+M\gamma+M+1}\cdot\frac{Nn}{3Nn+N\gamma+N+1}.
\end{align*}
Observe that this probability is given by $r_{2n}$ in \eqref{rr2}.
\end{itemize}

Consider now the case where $m$ is odd and write $m=2n+1,n\geq1$ (the simpler case $m=1$ will be treated separately). Remove all balls from urn A and place in urn A $Nn$ blue balls and $2Nn+N\gamma+2N+1$ red balls from the bath. In urn B we place $MNn+M-N$ blue balls and $N(2Mn+M\gamma+M+1)$ red balls. In urn R we place $M(n+1)$ blue balls and $2Mn+M\gamma+2M+1$ red balls. Draw one ball from urn A at random with the uniform distribution. If we get a blue ball then we go to urn B and draw a ball, while if we get a red ball then we go urn R and draw a ball. Again, we have three possibilities:
\begin{itemize}
\item If we get two blue balls in a row, i.e. one from urn A and then one from urn B, then we remove into the bath or add from the bath balls until we have $m-2=2n-1$ blue balls in urn A and start over. Therefore
$$
\mathbb{P}\left(X_1^{(1)}=2n-1 \, | \,  X_0^{(1)}=2n+1\right)=\frac{Nn}{3Nn+N\gamma+2N+1}\cdot\frac{MNn+M-N}{M(3Nn+N\gamma+N+1)}.
$$
Observe that this probability is given by $t_{2n+1}$ in \eqref{ts2}.
\item  If we get two red balls in a row, i.e. one from urn A and then one from urn R, then we remove into the bath or add from the bath balls until we have $m=2n+1$ blue balls in urn A and start over. Therefore
$$
\mathbb{P}\left(X_1^{(1)}=2n+1 \, | \,  X_0^{(1)}=2n+1\right)=\frac{2Nn+N\gamma+2N+1}{3Nn+N\gamma+2N+1}\cdot\frac{2Mn+M\gamma+2M+1}{3Mn+M\gamma+3M+1}.
$$
Observe that this probability is given by $s_{2n+1}$ in \eqref{ts2}.
\item If we get either a blue and a red ball or a red and a blue ball then we remove into the bath or add from the bath balls until we have $m-1=2n$ blue balls in urn A and start over. Therefore
\begin{align*}
\mathbb{P}\left(X_1^{(1)}=2n \, | \,  X_0^{(1)}=2n+1\right)&=\frac{Nn}{3Nn+N\gamma+2N+1}\cdot\frac{N(2Mn+M\gamma+M+1)}{M(3Nn+N\gamma+N+1)}\\
&\qquad\qquad+\frac{2Nn+N\gamma+2N+1}{3Nn+N\gamma+2N+1}\cdot\frac{M(n+1)}{3Mn+M\gamma+3M+1}.
\end{align*}
Observe that this probability is given by $r_{2n+1}$ in \eqref{rr2}.
\end{itemize}
For the case $m=1$ we only need to place in urn A $M$ blue balls and $M\gamma+2M+1$ red balls and draw one ball at random. If we get a blue ball (with probability $r_1$ in \eqref{rr2}) then we remove all balls from urn A and stop (since the state $0$ is absorbing). If we get a red ball (with probability $s_1$ in \eqref{ts2}) then we leave only $1$ blue ball in urn A and start over.

\bigskip

The urn model going with $P$ (on $\mathbb{Z}_{\geq0}$) is obtained by repeatedly alternating Experiments 1 and 2 in that order (see Figures \ref{fig1} and \ref{fig2} below, included here for the benefit of the reader). If we first perform Experiment 1 (for $m=2n$) we will end up with an urn with $2n-2, 2n-1$ or $2n$ blue balls. Now we perform Experiment 2 starting with $2n-2, 2n-1$ or $2n$ blue balls, in which case we may end with either $2n-1$ or $2n-2$ blue balls in the first case, $2n$ or $2n-1$ blue balls in the second case and $2n+1$ or $2n$ blue balls in the final third case. The situation is similar if we start with $m=2n+1$. We first perform Experiment 1 and we will end up with an urn with $2n-1, 2n$ or $2n+1$ blue balls. Now we perform Experiment 2 starting with $2n-1, 2n$ or $2n+1$ blue balls, in which case we may end with either $2n$ or $2n-1$ blue balls in the first case, $2n+1$ or $2n$ blue balls in the second case and $2n+2$ or $2n+1$ blue balls in the final third case. The combination of possibilities of these six cases yields the transition probabilities in \eqref{tpm} (see also \eqref{abcd}).

\begin{landscape}
\thispagestyle{empty}

\begin{figure}[h]
\begin{center}
\vspace{-0.0cm}
\includegraphics[height=13cm]{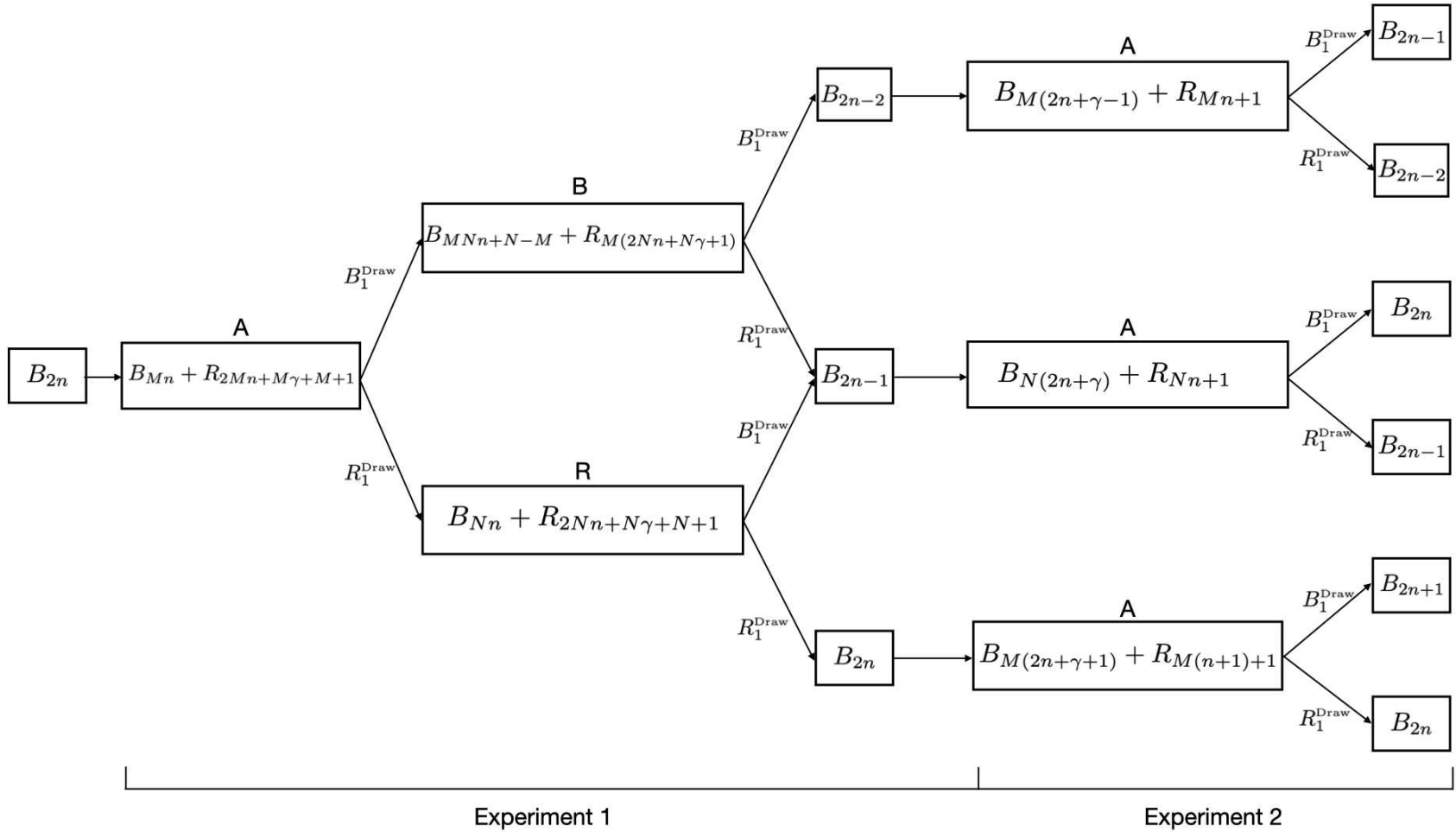}
\end{center}
\vspace{-0.2cm}
	\caption{A schematic of Experiments 1 and 2. The boxed regions represent the state or urn with $B_{\cdot}$ and $R_{\cdot}$ indicating the number of blue or red balls, respectively, contained within the urn. When a ball is drawn from an urn, this is indicated by $B_1^{\tiny{\mbox{Draw}}}$ or $R_1^{\tiny{\mbox{Draw}}}$. The initial state is $B_{2n}$.}
\label{fig1}
\end{figure}

\end{landscape}

\begin{landscape}
\thispagestyle{empty}

\begin{figure}[h]
\begin{center}
\vspace{-0.0cm}
\includegraphics[height=13cm]{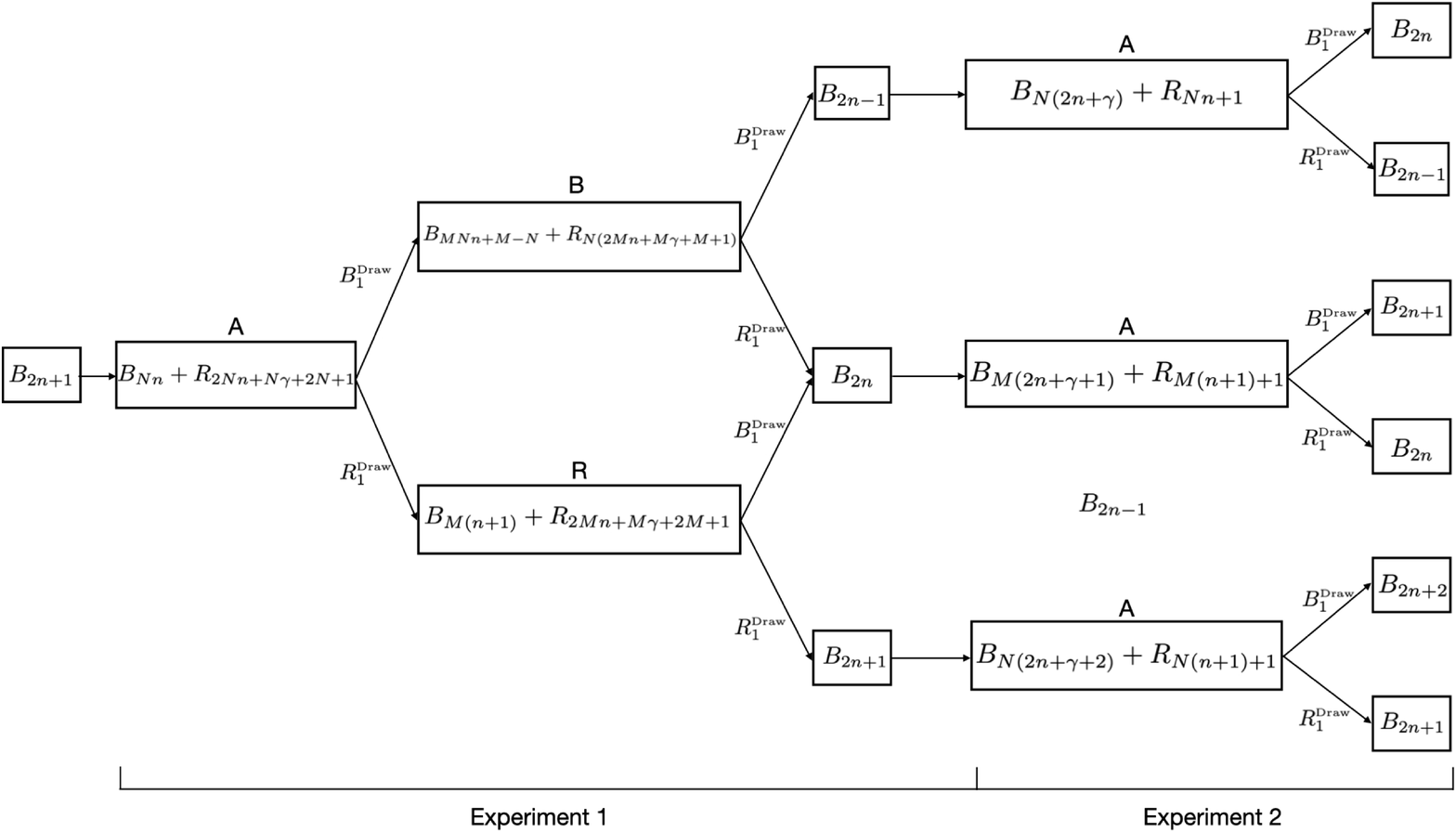}
\end{center}
\vspace{-0.2cm}
	\caption{A schematic of Experiments 1 and 2. The boxed regions represent the state or urn with $B_{\cdot}$ and $R_{\cdot}$ indicating the number of blue or red balls, respectively, contained within the urn. When a ball is drawn from an urn, this is indicated by $B_1^{\tiny{\mbox{Draw}}}$ or $R_1^{\tiny{\mbox{Draw}}}$. The initial state is $B_{2n+1}$.}
\label{fig2}
\end{figure}

\end{landscape}

\end{document}